\documentclass[graybox]{svmult}

\usepackage{type1cm}        % activate if the above 3 fonts are
\usepackage{makeidx}         % allows index generation
\usepackage{graphicx}        % standard LaTeX graphics tool
\usepackage{multicol}        % used for the two-column index
\usepackage[bottom]{footmisc}% places footnotes at page bottom

\usepackage{newtxtext}       % 
\usepackage[varvw]{newtxmath}       % selects Times Roman as basic font

\usepackage[colorlinks,citecolor=blue]{hyperref}

\usepackage{cite}

\usepackage{url}
\usepackage{graphicx}
\usepackage{float}
\usepackage{caption}
\usepackage{tikz}
\usepackage[ruled,vlined]{algorithm2e}
\usepackage{wrapfig}

\usetikzlibrary{calc,intersections,matrix,chains,positioning,decorations.pathreplacing,arrows}

\makeindex             % used for the subject index
\DeclareMathAlphabet{\bi}{OT1}{cmr}{bx}{it}

\newcommand{\R}{\mathbb{R}}

\newcommand{\x}{\bi x}

\newcommand{\bxi}{\boldsymbol{\xi}}
\renewcommand{\E}{{\bi E}}

\begin{document}

\title*{Finding the Shape of Lacunae of the Wave Equation Using Artificial Neural Networks}
\titlerunning{Artificial Neural Networks for the Wave Equation} %for an abbreviated version of
\author{Alina Chertock, Christopher Leonard, and Semyon Tsynkov}
\authorrunning{A. Chertock, C. Leonard, and S. Tsynkov} %for an abbreviated version of
\institute{Department of Mathematics, North Carolina State University, Raleigh, NC 27695, USA.
\email{chertock@math.ncsu.edu,cleonar@ncsu.edu,tsynkov@math.ncsu.edu} 
%\and Christopher Leonard \at Department of Mathematics, North Carolina State University, Raleigh, NC 27695, USA \email{cleonar@ncsu.edu}
%\and Semyon Tsynkov\at Department of Mathematics, North Carolina State University, Raleigh, NC 27695, USA \email{tsynkov@math.ncsu.edu}
}
\maketitle

\abstract{We apply a fully connected neural network to determine the shape of the lacunae in the solutions of the wave 
equation. Lacunae are the regions of quietness behind the trailing fronts of the propagating waves. The network is trained using a 
computer simulated data set containing a sufficiently large number of samples. The network is then shown to correctly reconstruct the 
shape of lacunae including the configurations when it is fully enclosed.}

\section{Introduction}\label{clt:sec1}

Consider the inhomogeneous scalar wave (d'Alembert) equation in the three-dimensional space (3D):
\begin{equation}
\label{clt:wave}
\frac{1}{c^2}\frac{\partial^2 u}{\partial t^2}-\Delta u=f(\x,t), \quad \x\in\R^3,\quad t\geqslant0,
\end{equation}
subject to zero initial conditions, and with the source term $f$ that is compactly supported on a bounded (3+1)D domain 
$Q_f\subset\R^3\times[0,+\infty)$.  The solution $u$ to equation \eqref{clt:wave} is given by the Kirchhoff integral:  
\begin{equation}
\label{clt:Kirchoff}
u(\x,t)=\frac{1}{4\pi}\iiint\limits_{|\x-\bxi|\leqslant ct}\frac{f(\bxi,t-|\x-\bxi|/c)}{|\x-\bxi|}d\bxi.
\end{equation}
The integration in \eqref{clt:Kirchoff} is performed in space over the ball of radius $ct$ centered at $\x$, but as $f$ is taken at 
retarded moments of time, this can be interpreted as integration in the (3+1)D space-time over the surface of a backward characteristic 
cone of equation (\ref{clt:wave}) (light cone of the past) with the vertex $(\x,t)$. This surface may or may not intersect with the 
support $Q_f$ of the right-hand side $f$. If there is no intersection, then  $u(\x,t)=0$, which implies, in particular, that the solution 
$u=u(\x,t)$ of equation (\ref{clt:wave}) will have a lacuna (secondary lacuna in the sense of Petrowsky \cite{petrowsky-45}): % that we denote $\Lambda$:
\begin{equation}
\label{clt:lac}
u(\x,t)\equiv0\ \ \
\forall(\x,t)\in\!\!
\bigcap_{(\bxi,\tau)\in Q_f}\!
\left\{\bigl.(\tilde \x,\tilde t)\bigr||\tilde \x-\bxi|<c(\tilde t-\tau),\
\tilde t>\tau\right\} \stackrel{\text{def}}{=}\Lambda.
\end{equation}
Mathematically, %formula (\ref{clt:lac}) implies that 
the lacuna $\Lambda$ is the intersection of all forward characteristic cones (i.e., light cones of the future) of the wave equation 
(\ref{clt:wave}) once the vertex of the cone sweeps the support $Q_f$ of the  right-hand side $f(\x,t)$. From the standpoint of  physics, 
$\Lambda$ is the part of space-time where the waves generated by a compactly supported source have already passed and the solution has 
become zero again. The primary lacuna (as opposed to  secondary lacuna (\ref{clt:lac})) is the part of  
space-time ahead of the propagating fronts where the waves have not reached yet. % Hereafter, we will be focusing on the secondary 
%lacunae, and  refer to those as merely  lacunae.

The phenomenon of lacunae is inherently three-dimensional  (more precisely, it pertains to  spaces of odd dimension). The surface of the lacuna includes 
the trajectory of aft (trailing) fronts of the propagating waves. The existence of {sharp aft fronts} in odd-dimension spaces is known as 
the (strong) Huygens' principle, as opposed to the so-called wave diffusion, which takes place in spaces of even dimension 
\cite{vlad,courant2}.

The question of identifying the hyperbolic equations and systems that admit the diffusionless propagation of waves has been first  
formulated by Hadamard \cite{hadamard1,hadamard2,hadamard3}. He, however, did not know any other examples besides the d'Alembert equation 
(\ref{clt:wave}). The notion of lacunae was introduced and studied by Petrowsky in \cite{petrowsky-45}, where
conditions for the coefficients of hyperbolic equations that guaranteed the existence of lacunae have been obtained (see also \cite[Chapter VI]{courant2}).
Subsequent developments can be found in \cite{abg1,abg2}. However, since work \cite{petrowsky-45} no other constructive examples of either
scalar equations or systems that satisfy the Huygens' principle have been found except for the wave equation (\ref{clt:wave}) and its 
equivalents. Specifically, it was shown in \cite{matthisson} that in the standard $(3+1)$D space-time with Minkowski metric, the only 
scalar hyperbolic equation that has lacunae is the wave equation (\ref{clt:wave}). The first examples of nontrivial  diffusionless
equations (i.e., irreducible to the wave equation) were constructed in \cite{stell1,stell2,stell3}, but the space must be $\R^d$ for odd 
$d\geqslant5$. Examples of nontrivial diffusionless  systems (as opposed to scalar equations) in the standard Minkowski $(3+1)$D 
space-time were presented in \cite{schimming,belger,gunther}, as well as examples of nontrivial scalar Huygens' equations in a $(3+1)$D 
space-time  equipped with a different  metric (the  plane wave metric that contains off-diagonal terms), see
\cite{belger,gunther,gunther3}. It was shown in \cite{lax-78} that the wave equation on the $d$-dimensional sphere, where $d\geqslant3$ 
is odd, satisfies the Huygens' principle; this spherical wave equation can be transformed to the Euclidean wave equation locally, but not 
globally.

While the examples of nontrivial diffusionless equations/systems built in \cite{stell1,stell2,stell3,schimming,belger,gunther,gunther3} 
are  primarily of a theoretical interest, the original wave equation (\ref{clt:wave}) accounts for a variety of physically relevant (albeit 
sometimes simplified) models in acoustics, electromagnetism, elastodynamics, etc. Accordingly, understanding the shape of the lacunae 
(\ref{clt:lac}) is of interest for the aforementioned application areas as lacunae represent the regions of ``quietness'' where the 
corresponding wave field is zero. 

The objective of our work is to determine the shape of the lacunae in the solutions of the wave equation using the fully connected 
artificial neural networks. In the current paper, we adopt a simplified scenario to construct, test, and verify the proposed machine 
learning approach. Specifically, while the true phenomenon of lacunae is 3D and applies to solutions given by the Kirchhoff integral 
(\ref{clt:Kirchoff}), hereafter we conduct the analysis and simulations in a one-dimensional (1D) setting. The domain of dependence for 
$u(\x,t)$ determined by the Kirchhoff integral is the surface of the backward light cone. To mimic that in 1D, we consider the function 
$u=u(x,t)$
and define its domain of dependence as the sum of two backward propagating rays: 
\begin{equation}
\label{clt:dd1}
\{(\xi,\tau): \xi-x=\pm c(\tau-t), \ \ \xi\in\R, \ \ \tau\leqslant t\}\stackrel{\text{def}}{=}\mathcal{L}(x,t).
\end{equation}
We do not need a full specification of $u$ for our subsequent considerations. We only need a sufficient condition for $u(x,t)$ to be equal to zero, which we take as % $u(x,t)=0$ if $\mathcal{L}(x,t)\cap Q_f=\emptyset$,
$$ u(x,t)=0 \quad\text{if}\quad \mathcal{L}(x,t)\cap Q_f=\emptyset,$$
where $Q_f$ is a given bounded domain in (1+1)D space-time: $Q_f\subset \R\times[0,+\infty)$.
Accordingly,  the function $u$ is going 
to have a lacuna:
\begin{equation}
\label{clt:lac1}
u(x,t)\equiv0\ \ \
\forall(x,t)\in\left\{\bigl.(\tilde x,\tilde t):\mathcal{L}(\tilde x,\tilde t)\cap Q_f=\emptyset\right\} \stackrel{\text{def}}{=}\Lambda_1(Q_f).
\end{equation}
The lacuna $\Lambda_1(Q_f)$ combines both the secondary and primary lacuna as per the discussion in Section~\ref{clt:sec1}. A purely secondary lacuna would be given by [cf.\  (\ref{clt:lac})]
\begin{equation}
\label{clt:lac1-2}
u(x,t)\equiv0\ \ \
\forall(x,t)\in\!\!
\bigcap_{(\xi,\tau)\in Q_f}\!
\left\{\bigl.(\tilde x,\tilde t):|\tilde x-\xi|<c(\tilde t-\tau),\
\tilde t>\tau\right\} \subset \Lambda_1(Q_f).
\end{equation}

We emphasize that the proposed 1D construct is not accurate on the substance. Its only purpose is to provide an inexpensive testing 
framework for the neural networks  described in the paper. This construct is designed as a direct counterpart of the physical 
3D setting and does not represent a true solution of the 1D wave equation:
\begin{equation}
\frac{1}{c^2}u_{tt}-u_{xx}=f(x,t),\quad x\in\mathbb{R},\ t>0.
\label{clt:1D}
\end{equation} 
 The solution of (\ref{clt:1D}) subject to zero initial conditions is given by the d'Alembert integral [cf.\  the
Kirchhoff integral (\ref{clt:Kirchoff})]:
\begin{equation}
\label{clt:d'Alembert}
u(x,t)=\frac{c}{2}\int\nolimits_0^td\tau\int\nolimits_{x+c(\tau-t)}^{x-c(\tau-t)}f(\xi,\tau)d\xi.
\end{equation}
Unlike  (\ref{clt:dd1}), the domain of dependence for the 1D solution (\ref{clt:d'Alembert}) contains not only the two rays, but 
the entire in-between region as well: 
\begin{equation*}
%\label{eq:dd2}
\{(\xi,\tau):\ x + c(\tau-t) \leqslant \xi\leqslant x - c(\tau-t), \ \ \tau\leqslant t\} \stackrel{\text{def}}{=}\mathcal{D}(x,y).
\end{equation*}
Therefore, the solution $u=u(x,t)$ defined by (\ref{clt:d'Alembert}) will not, generally speaking, have  secondary lacunae as presented by 
(\ref{clt:lac1-2}).\footnote{The  1D case is special in the sense that while the  dimension of the space is odd, the wave 
equation (\ref{clt:1D}) is not Huygens' when driven by a source term. It may, however, demonstrate a Huygens' behavior with respect to the initial data 
\cite{vlad}.
}

The paper is organized as follows. In Section \ref{clt:sec2}, we provide 
a description of the numerical algorithm to build the data set. In Section~\ref{clt:sec3}, we introduce the neural networks (NN) and and discuss their training. In Section~\ref{clt:sec4}, we present the numerical results. Section~\ref{clt:sec5} contains our conclusions and identifies directions for future work.

\section{Construction of the Data Set}\label{clt:sec2}
In this section, we describe the construction of the data set needed for the network training. To this end, we introduce a computational 
domain $\Omega:=[a,b]\times[0,T]$ and discretize it in space and time using a uniform spatial $x_1,\ldots,x_{N_x}$ and temporal 
$t_1,\ldots,t_{N_t}$ mesh so that
\begin{align*}
x_j&=a+(j-1)\Delta x,\quad j=1,\ldots,N_x\\
t_n&=(n-1)\Delta t,\quad n=1,\ldots,N_t,
\end{align*}
where $\Delta x=\frac{b-a}{N_x-1}$ and $\Delta t=\frac{T}{N_t-1}$. 

%\begin{figure}[H]
\begin{wrapfigure}[14]{l}{.40\textwidth}
%\centering
\vspace{-.4cm}
\hspace{-3mm}
\includegraphics[trim= 50 25 0 40,clip,width=0.45\textwidth]{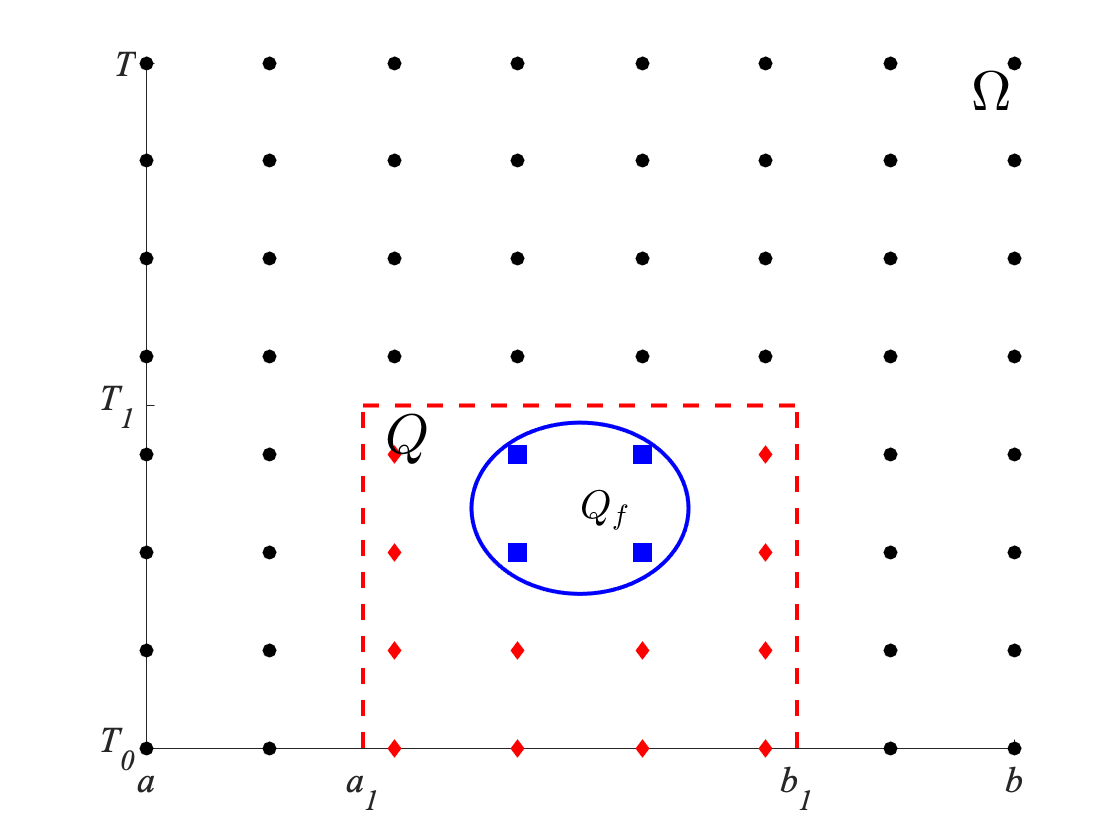}
\caption{The computational domain $\Omega$ with  subdomain $Q$ inside the dotted red line, and $Q_f$ inside the blue circle.}
\label{clt:fig:domain}  
%\end{figure}
\end{wrapfigure} 
We assume that the domain $Q_f$ that defines the lacuna $\Lambda_1$ in (\ref{clt:lac1}) lies inside a subdomain $Q\subset\Omega$, which,
for the simplicity of implementation, is taken as a box $Q=[a_1,b_1]\times[T_0,T_1]\subseteq\Omega$, where $a\leq a_1<b_1\leq b$ and 
$0\leq T_0< T_1\leq T$. We can then denote the nodes inside  $Q$ by
\begin{align*}
x_{j_\ell}&=x_{j_1+(\ell-1)}, \quad \ell=1,\ldots,N_x^{'}\leq N_x,\\ 
t_{n_p}&=t_{n_1+(p-1)}, \quad p=1,\ldots,N_t^{'}\leq N_t,  
\end{align*}
where $x_{j_1}$ and $t_{n_1}$ are the smallest $x_j$ and $t_n$ values that are inside the the set $Q$ and $x_{j_1+N^{'}_x-1}$ and 
$t_{n_1+N^{'}_t-1}$ are the largest $x_j$ and $t_n$ values that are inside the the set $Q$, see Figure~\ref{clt:fig:domain}.

\begin{figure}[ht!]
\begin{center}
	\includegraphics[trim= 50 25 0 50,clip,width=0.5\textwidth]{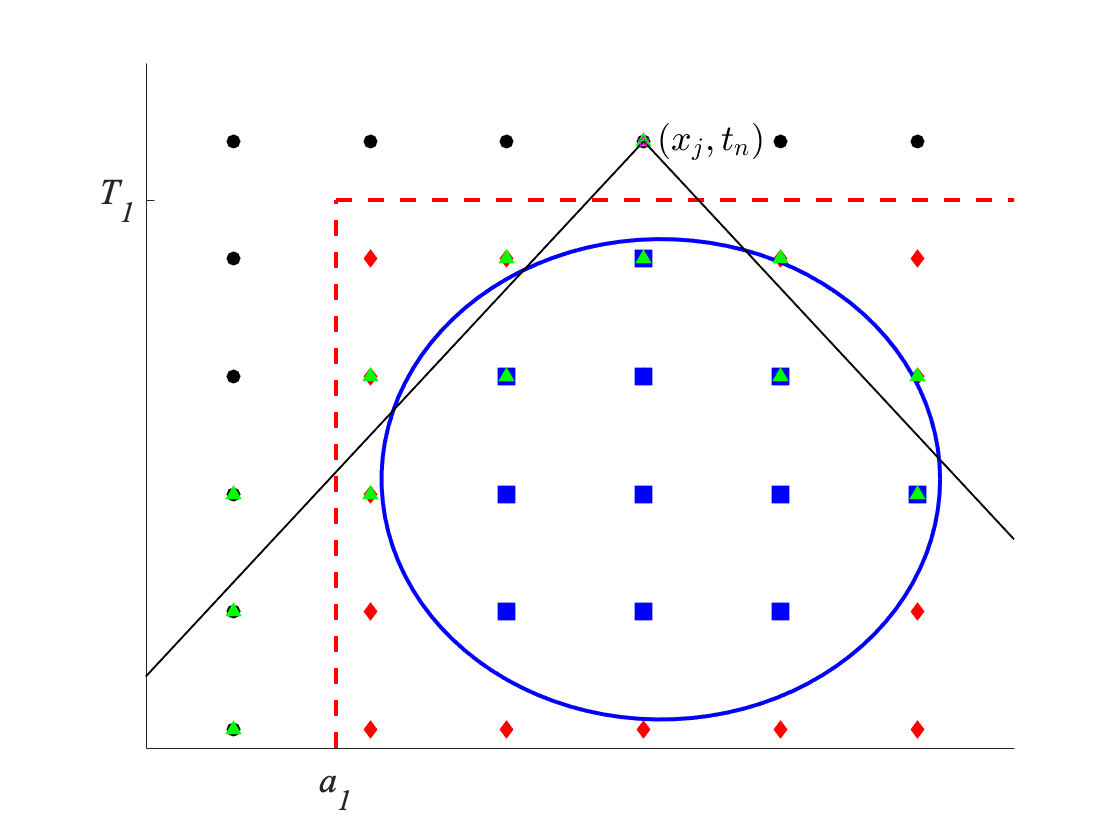}
	\caption{Characteristic lines from the point $(x_j,t_n)$. Green triangles indicate nodes within $\Delta x$ of the characteristic lines at each time step. If there exist a green triangle and blue square at the same node, then $\Psi_{j,n}=1$, else $\Psi_{j,n}=-1$.} 
	\label{clt:fig:nodes}
	\end{center}
\end{figure}

The shape of the lacuna can then be determined by identifying the set of nodes, for 
which the characteristic lines $\mathcal{L}(x_j,t_n)$ emerging from the node 
$(x_j,t_n)\in\Omega,\ j=1,\ldots,N_x,\ n=1,\ldots,N_t$ (see formula \eqref{clt:dd1}), pass through the domain $Q_{f}$.  
We therefore construct $M$ training data sets by implementing the following algorithm:

\newpage 
\begin{programcode}{Algorithm 1}
\begin{enumerate}
\item
{\it Start}: Introduce the computational domain $\Omega$ and its discretization by a uniform mesh  
$(x_j,t_n)\in\Omega,\ j=1,\ldots,N_x,\ n=1,\ldots,N_t$, and identify the set of nodes 
$(x_{j_\ell},t_{n_p})\in Q,\ \ell=1,\ldots,N^{'}_x,\ p=1,\ldots,N^{'}_t$.
\item
{\it Iterate}: For $m=1,2, \ldots,M$
\begin{enumerate}
\item
Generate a random positive integer $I^{(m)}$ and a set of domains $Q_{f^{(m)}_i}\subset Q$ 
%that represent supports of some functions 
%$f^{(m)}_i(x,t)$ 
for each $i=1,\ldots,I^{(m)}$, and define
\begin{equation}\label{clt:Qf}
Q_{f^{(m)}}=\Big(\bigcup\limits_{i=1}^{I^{(m)}} Q_{f^{(m)}_i}\Big).
\end{equation}
\item
Construct the following matrices $\Phi^{(m)}$ and $\Psi^{(m)}$:
\begin{itemize}
\item
Matrix $\Phi^{(m)}$ with entries $\phi_{\ell,p}^{(m)}$ that indicate for each node $(x_{j_\ell},t_{n_p})\in Q$ whether or not it belongs 
to the domain $Q_{f^{(m)}}$, namely,
\begin{equation*}
\phi_{\ell,p}^{(m)}=\left\{\begin{aligned}
&1,&&(x_{j_\ell},t_{n_p})\in Q_{f^{(m)}},\\
&-1,&&\mbox{otherwise.}
\end{aligned}\right.
%\label{Bjn}
\end{equation*} 
\item 
Matrix $\Psi^{(m)}$ with entries $\psi_{j,n}^{(m)}$ indicates whether or not the characteristic lines $\mathcal{L}(x_j,t_n)$ intersect 
with the domain $Q_{f^{(m)}}$,  namely,
\begin{equation}
\psi_{j,n}^{(m)}=\left\{\begin{aligned}
&-1,&& (x_j,t_n)\in\Lambda_1(Q_{f^{(m)}}),\\
&1,&&\mbox{otherwise.}
\end{aligned}\right.
\label{clt:Psi_jn}
\end{equation}
In practice, we check whether there is a node $(x_{j_\ell},t_{n_p})\in Q_{f^{(m)}}$ and point $(\xi,t_{n_p})\in{\mathcal{L}(x_j,t_n)}$, such 
that $|\xi-x_{j_\ell}|<\Delta x$, in which case $\psi_{j,n}^{(m)}=1$; otherwise $\psi_{j,n}^{(m)}=-1$. Note that, the point $(\xi,t_{n_p})$ is not, generally speaking,  a grid node. The constriction of $\Psi^{(m)}$ is visualized in Figure~\ref{clt:fig:nodes}. 
\end{itemize}
\end{enumerate}
\item
{\it Form data set:} Store the set $\{\Phi^{(m)},\Psi^{(m)}\}_{m=1}^M$.
\end{enumerate}
\end{programcode}

%\runinhead{Remark}
%In one spacial dimension, this approach to constructing the matrices $\Psi^{(m)},\ m=1,\ldots,M,$ could be done with a more exact approach by 
%finding the lowest and highest characteristic lines coming from the the connected regions $Q_{f^{(m)}_i}$, and finding the nodes between 
%those lines. However, we believe a similar approach to the one described in algorithm 1 will work better in 3D so we use it here as 
%well.

\section{Construction of the Neural Network and Training}\label{clt:sec3}
%To train a fully connected feed forward neural network $N_\Theta:\Z_2^{N^{'}_x\times N^{'}_t}\to \R^{N_x\times N_t}$, we solve the 
%following optimization problem:
%\begin{equation}\label{loss}
%\begin{split}
%\Theta^* &=\underset{\Theta}{\arg\min}\; L(\Theta)\\
% L(\Theta) &= \frac{1}{M}\sum_{m=1}^M\|\Psi^{(m)}-N_\Theta(\Phi^{(m)})\|,
% \end{split}
%\end{equation}
%where $\|\cdot\|$ is a given norm.

\begin{figure}[ht!]
  \begin{center}
    \begin{tikzpicture}
        [greennode/.style={rectangle, draw=green!60, fill=green!5, very thick, minimum width=0.6cm, minimum height = 2cm,scale=1},
        bluenode/.style={rectangle, draw=blue!60, fill=blue!5, very thick, minimum width=.6cm, minimum height = 2cm,scale=1},
        yellownode/.style={rectangle, draw=yellow!60, fill=yellow!5, very thick, minimum width=0.6cm, minimum height = 2cm,scale=1},
        graynode/.style={rectangle, draw=gray!60, fill=gray!5, very thick, minimum width=0.6cm, minimum height = 2cm,scale=1}]
        \node at (1.0,5) [greennode]        (i1)    {\rotatebox{90}{input}};
        %\node[graynode]         (r1)       [right=0.8cm of i1]{\rotatebox{90}{reshape ($M_xM_t$)}};
        \node[bluenode]         (h1)       [right=1.2cm of i1]{\rotatebox{90}{$\sigma_1$}};
        \node[bluenode]         (h2)       [right=1.2cm of h1]{\rotatebox{90}{$\sigma_2$}};
        \node[]       (o1)       [right=.3cm of h2]{{$\dots$}}; 
        \node[bluenode]        (h3)       [right=1.2cm of h2]{\rotatebox{90}{$\sigma_k$}};
        \node[yellownode]       (o1)       [right=1.2cm of h3]{\rotatebox{90}{$\sigma_{K+1}$ (output)}};
        %\node[graynode]         (r2)       [right=0.8cm of o1]{\rotatebox{90}{reshape ($N_x\times N_t$)}};
        
        %\draw[->] (i1.east) -- (r1.west)
        node[pos=0.5,sloped,above] {};
        \draw[->] (i1.east) -- (h1.west)
        node[pos=0.5,sloped,above] {$A_1$};
        \draw[->] (h1.east) -- (h2.west)
        node[pos=0.5,sloped,above] {$A_2$};
        
        \draw[->] (h3.east) -- (o1.west)
        node[pos=0.5,sloped,above] {$A_{K+1}$};
        %\draw[->] (o1.east) -- (r2.west)
        node[pos=0.5,sloped,above] {};
        
        \useasboundingbox (0,5) rectangle (10,5);
    \end{tikzpicture}\\
    \caption{Neural Network Architecture}
    \label{clt:arch}
  \end{center}
\end{figure}
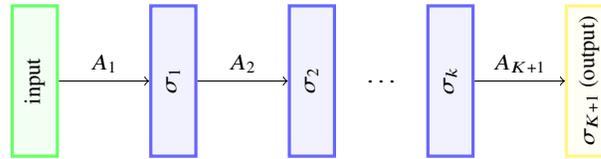

In this paper, we use a fully connected feed forward neural network $N_\Theta:\R^{N^{'}_x\times N^{'}_t}\to \R^{N_x\times N_t}$ to approximate the shape of the lacunae. To train $N_\Theta$, we search for the parameter set $\Theta$ that minimizes the loss function:
\begin{equation}\label{clt:loss}
 L(\Theta) = \frac{1}{M}\sum_{m=1}^M\|\Psi^{(m)}-N_\Theta(\Phi^{(m)})\|,
\end{equation}
where $\|\cdot\|$ is a given norm.
The fully connected neural network with $K$ hidden layers is the composition of functions
\begin{equation}
\label{clt:NN}
N_\Theta = N_{\Theta_{K+1}}\circ N_{\Theta_{K}}\circ...\circ N_{\Theta_1}, 
\end{equation}
with    
\begin{equation*}%\label{NN}
N_{\Theta_k}(z) = \sigma_k(A_k),\quad A_k=W_kz+b_k,\quad k=1,2,...,K+1,
\end{equation*}
where $A_k$ is the connection between layers $k-1$ and $k$, and $\sigma_k$ is the activation for the layer $k$. Here, $k=0$ corresponds to 
the input layer and $k=K+1$ corresponds to the output layer, and we assume that there are $w_k$ nodes on each layer $k=0,\ldots,K+1$. For
each connection between layers,  $W_k\in\R^{w_k\times w_{k-1}}$ is the weight matrix and $b_k\in\R^{w_k}$ is the bias vector. The 
parameters in the sets $\Theta_k$ are the element values for the weight matrix $W_k$ and bias vector $b_k$, and 
$\Theta=\cup_{k=1}^{K+1}\Theta_k$ is the parameter set that we train the neural network to learn by minimizing the loss function \eqref{clt:loss}. 
%
%Note that the input and output of this neural network is a vector so one needs to reshape the matrices $\Phi^{(m)}$ and $\Psi^{(m)}$, $m=1,2,...,M$, into vectors of size $N^{'}_xN^{'}_t=w_0$ and $N_xN_t=w_{K+1}$ respectively.    
%
We can represent the architecture of a fully connected feed forward neural network as  shown in Figure \ref{clt:arch}.

For each $k=1,2,...,K+1$, the activation function $\sigma_k$ is usually a nonlinear function that is applied to each element of its input vector. There are many activation functions that can be used, some of the most popular activation functions used include the rectified linear unit (ReLU), sigmoid functions, softmax, etc \cite{NEURIPS2019_9015}.

We train the neural network $N_\Theta$ with the following algorithm:
\begin{programcode}{Algorithm2}
\begin{enumerate}
\item
{\it Start}: Randomly split the data set into two sets, the training set of size $M_{tr}$, $\{\Phi^{(m_i)},\Psi^{(m_i)}\}$, $i=1,2,\ldots,M_{tr}$ and the validation set of size $M_{val}$, $\{\Phi^{(m_l)},\Psi^{(m_l)}\}$, $l=1,2,\ldots,M_{val}$. 
\begin{enumerate}
\item The training set is used in the optimization algorithm to find the parameter set $\Theta$ that minimizes \eqref{clt:loss}. 
\item The validation set is used to check that the neural network can generalize to data that are not in the training set.
\end{enumerate}

\item
{\it Set hyperparameters}: loss function, optimizer, learning rate, number of epochs, batch size, number of hidden layers, number of nodes per hidden layer, and activation functions.
\begin{enumerate}
\item loss function: The function that is to be minimized given the training data set, see \eqref{clt:loss}. 
\item optimizer: The algorithm that is used to try and find the global minimum of the loss function. 
\item learning rate: The initial step size that the optimization algorithm takes. Depending on the optimization algorithm, the step sizes may change between steps.  
\item number of epochs: The number of times the optimization algorithm goes through the entire training data set.
\item batch size: The number of samples from the training data that will propagate through the network for each update of the parameters. 
\item number of hidden layers: $K$.
\item number of nodes per each hidden layer: $w_k,\ k=1,\ldots,K$. 
\item activation functions: $\sigma_k,\ k=1,\ldots,K+1$. 
\end{enumerate}

%\item
%Set $best\_loss=\infty$.
\item Randomly initialize the parameter set: $\Theta^{(0)}$. 
\item
{\it Iterate}: For $e=1,2, \ldots, \text{number of epochs}$
\begin{enumerate}
\item
Randomly split the training set into $\beta$ separate batches  
\item
{\it Iterate}: For batch=1,2,$\ldots$,$\beta$

\begin{itemize}
\item 
Update the neural network parameters using one step of the optimization algorithm. 
\end{itemize}
\item 
With the current parameter set $\Theta^{(e)}$, evaluate the output of $N_{\Theta^{(e)}}$ (see formula (\ref{clt:NN})) applied to all of the input values $\Phi^{(m_l)}$, $l=1,2,...,M_{val}$, in the validation set,  and calculate the loss value $L(\Theta^{(e)})$ defined in \eqref{clt:loss}. If this loss value is smaller than at the end of every other epoch before it, set $\Theta=\Theta^{(e)}$.
%\item
%Set $Loss_e$ to be the loss function applied to the validation set. 
%\item
%If $Loss_e<best\_loss$ then save parameter set $\Theta$ and set $best\_loss=Loss_e$
 
\end{enumerate}
\item
{\it Retrieve parameters}: Return $N_{\Theta}$ with the parameter set $\Theta$ obtained in step 3c.
\end{enumerate}
%\end{verbatim}
\end{programcode}

\runinhead{Remark} To apply the neural network to any input matrix $\Phi$, one needs to reshape the matrix into a vector of size $N^{'}_xN^{'}_t$. One can then reshape the vector output to be a matrix of the size of $\Psi$, which is $N_x\times N_t$.  

\runinhead{Remark} The training set is usually about 80\% of the total data set and the other 20\% is the validation set.

\runinhead{Remark} Finding good hyperparameters is often done by a running multiple trials of the training algorithm using different hyperparameters and identifying which one results in the best outcome.  The specific hyperparameters we have used are presented in Section~\ref{clt:sec4}.

\section{Numerical Results}\label{clt:sec4}
In this section, we conduct a number of numerical experiments to demonstrate the performance of the machine learning approach to detect the lacunae. All of our models are built and trained using the open source library PyTorch \cite{NEURIPS2019_9015}.

To define the set $Q_{f^{(m)}_i}$ in \eqref{clt:Qf}, we draw uniformly distributed values $\x_i^{(m)}\sim\mathcal{U}(a_1,b_1)$, $t_i^{(m)}\sim\mathcal{U}(0,T_1)$, and $r_i^{(m)}\sim\mathcal{U}(0,R)$ for $i=1,2,..I^{(m)}$, and then let 
\begin{equation}\label{clt:qfi}
Q_{f^{(m)}_i}=\{(\x,t)\in[a_1,b_1]\times[0,T_1]: (\x-\x_i^{(m)})^2+(t-t_i^{(m)})^2\leq (r_i^{(m)})^2 \}.
\end{equation}
Recall that $I^{(m)}$ is the number of sets on the right hand side of \eqref{clt:Qf}.
%\begin{equation}\label{f}
%f^{(m)}=\sum_{i=1}^{I^{(m)}} f^{(m)}_i, \quad i=1,\ldots,M.
%\end{equation} 
As mentioned in Algorithm 1, this integer is randomly generated, and in our numerical examples it is an integer between 1 and 4. If $I^{(m)}>1$, then we may have a disconnected set $Q_{f^{(m)}}$. 

%Equipped with $Q_f^{(m)}$ where $Q_{f_i}$ is defined above, we can create the matrix $\Phi$ as described in algorithm1 and train the neural network $N_\Theta$ such that    
%$$
%N_\Theta(\Phi)\approx \Psi
%$$
%where $\Psi$ is the matrix representing the location of the lacunae $\Lambda_1(Q_f)$ and its compliment $(\Lambda_1(Q_f))^c$. 

For all of the examples, the computational domains are $\Omega=[-20,20]\times[0,20]$ and $Q=[-10,10]\times[0,10]$. We discretize $\Omega$ and $Q$ such that $N_x=64$, $N_t=64$, $N_x^{'}=32$, and $N_t^{'}=32$. The max radius in \eqref{clt:qfi} we use is $R=5$. We generated $M=10{,}000$ data pairs for our training process, randomly splitting the data such that $M_{tr}=8{,}000$ and $M_{val}=2{,}000$. For the hyperparameters in Algorithm 2 we use: 
\begin{itemize}
	\item[a.] $L(\Theta) = \frac{1}{M}\sum_{m=1}^M\|\Psi^{(m)}-N_\Theta(\Phi^{(m)})\|_{F}$, where $\|\cdot\|_F$ is the Frobenius norm. 
	\item[b.] The Adam optimizer \cite{KB14}. 
	\item[c.] Learning rate $=10^{-4}$.
	\item[d.] Number of epochs $= 200$.
	\item[e.] Batch size $=32$.
	\item[f.] $K = 3$.
	\item[g.] $w_k=256$ for all $k=1,\ldots,K$.
	\item[h.] $\sigma_k=LeakyReLU$ \cite{NEURIPS2019_9015} for $k=1,\ldots,K$ and $\sigma_{K+1}=\tanh$.
\end{itemize}   

Equipped with all the parameters and data, we train a neural network $N_\Theta$ such that for a given set 
$Q_f$ represented by the matrix $\Phi$ as described in Algorithm 1, the neural network will produce    
$$
N_\Theta(\Phi)\approx \Psi,
$$
where $\Psi$ is the matrix representing the location of the lacunae $\Lambda_1(Q_f)$ and its compliment $(\Lambda_1(Q_f))^c$.

\subsection{Example 1: Case $I^{(m)}=1$}
In the first example, we consider a particular case where $I^{(m)}=1$ for each $m=1,\ldots,M$, and  $M=10{,}000$. After training the neural network, $N_\Theta$, we apply it to a test set that is independent of the training and validation sets. For a given input matrix $\Phi$, let $\Psi^{ref}$ denote the reference solution matrix with elements $\psi^{ref}_{j,n}\in\{-1,1\}$, $j=1,\ldots,N_x$, $n=1,\ldots,N_t$, 
determining whether or not the node $(x_j,t_n)$ is in the lacuna, see \eqref{clt:Psi_jn}. Then for the neural network solution, $\Psi^{NN}\coloneqq N_\Theta(\Phi)\in(-1,1)$, we state that it correctly identifies that the node $(x_j,t_n)$ is in the right set, $\Lambda_1(Q_f)$ or $\Lambda_1(Q_f)^c$, if 
\begin{equation*}%\label{correct}
	(\psi^{NN}_{j,n} \leq 0 \text{ and } \psi^{ref}_{j,n}=-1)\quad
	\text{or}\quad  (\psi^{nn}_{j,n} > 0 \text{ and } \psi^{ref}_{j,n}=1),
\end{equation*}    
and incorrectly identifies which set the node is in if 
\begin{equation*}\label{incorrect}
	(\psi^{NN}_{j,n} > 0 \text{ and } \psi^{ref}_{j,n}=-1)\quad
	\text{or}\quad (\psi^{nn}_{j,n} \leq 0 \text{ and } \psi^{ref}_{j,n}=1).
\end{equation*} 
We can then determine how accurate the neural network is by calculating 
\begin{equation}\label{clt:acc}
	accuracy = \frac{\text{\# of nodes correctly identified}}{\text{total \# of nodes}}.
\end{equation}

\begin{figure}[H]
    %\hspace{-1.2cm} %w=1.2
\hspace{-0.0cm} %w=1.6
%trim= 100 0 0 40,clip,
\includegraphics[trim=80 0 0 0,clip,height=.86\textwidth,width=1.08\textwidth]{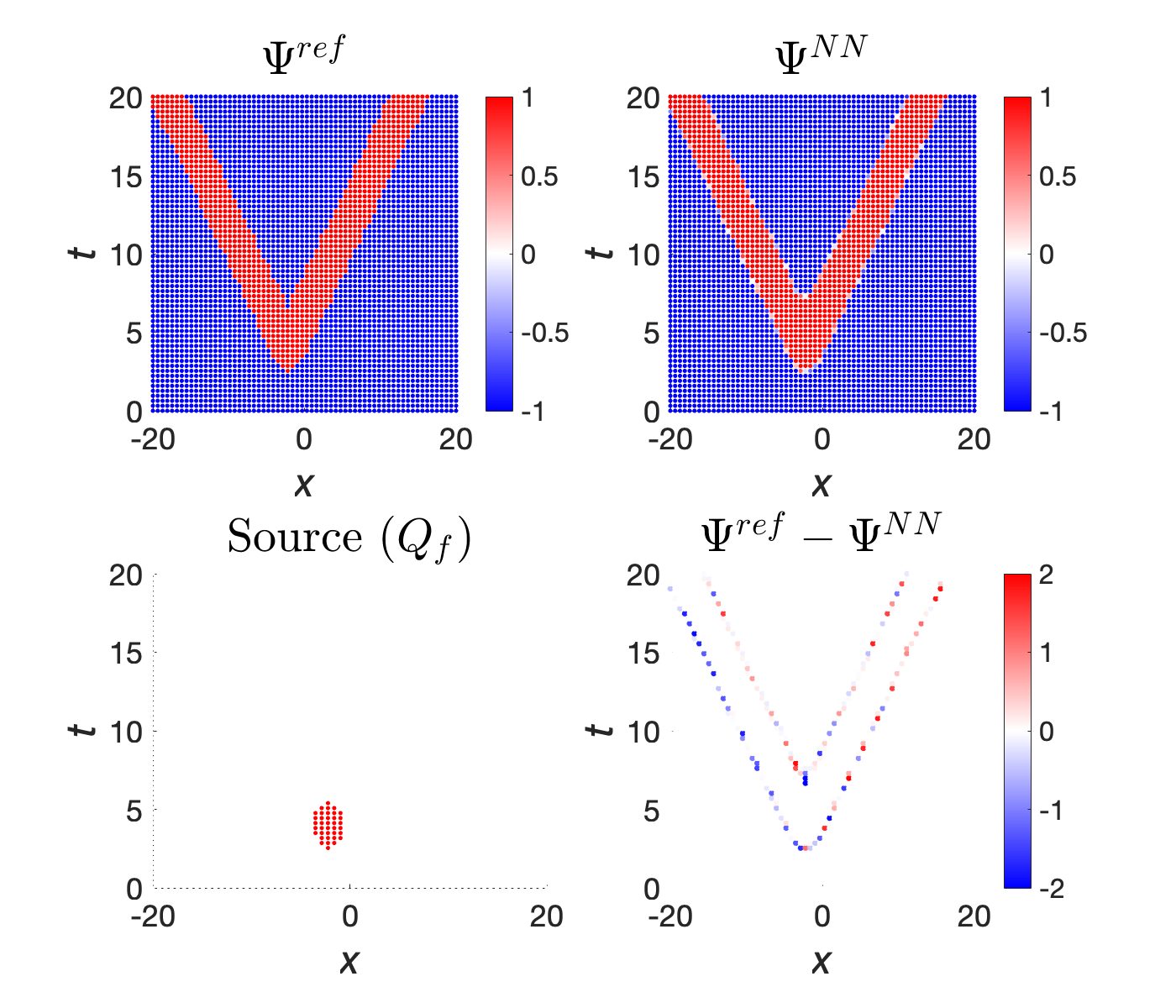}
    \caption{Reconstruction of the shape of the lacuna (\ref{clt:lac1}) by neural network (\ref{clt:NN}) in the case $I^{(m)}=1$. Top left: Reference solutions $\Psi^{ref}$.  Top right: Neural Network solutions $\Psi^{NN}$. Bottom left: The sets $Q_f$. Bottom right: $\Psi^{ref}-\Psi^{NN}$.
\label{clt:fig:single}}
\end{figure}

From $1{,}000$ test cases, the neural network was able to identify which set each node belonged to with approximately $99.12\%$ accuracy. Figure \ref{clt:fig:single} shows an examples taken from the test set. The top left represents the values for the reference solutions $\Psi^{ref}$, the top right represents the values of the neural network solution $\Psi^{NN}$, the bottom left shows nodes in the set $Q_f$, and the bottom right represents the difference $\Psi^{ref}-\Psi^{NN}$. For the top graphs, the values range from $[-1,1]$ with $-1$ indicating nodes in $\Lambda_1(Q_f)$ with blue dots and $1$ indicating nodes in $\Lambda_1(Q_f)^c$ with red dots. For the reference solution, the element values are either -1 or 1 so all the nodes are either dark blue or dark red respectively. The neural network has values between $-1$ and $1$, thus the node colors in its graph might be different shades of red, blue, and white. On the plot for $\Psi^{ref}-\Psi^{NN}$, the white nodes indicates that the two solution are close to each other, the red nodes indicates that the $\Psi^{ref}$ is greater than $\Psi^{NN}$ and the blue nodes indicate $\Psi^{ref}$ is less than $\Psi^{NN}$. Note that, the interior of the sets $\Lambda_1(Q_f)$ and $\Lambda_1(Q_f)^c$ for the neural network solution are clearly defined as very close to -1 or 1, but there is some uncertainty from the neural networks along the boundary. It is expected that the neural network would in general have more difficulty learning the edges of these sets. 
      
\vspace{-1.0cm}      
\subsection{Example 2:  Case $1\leq I^{(m)}\leq4$} 
\begin{figure}[H]
	%width=1.0\textwidth
    %\hspace{0cm}  % w=1.6
    \includegraphics[trim=80 0 0 0,clip,height=0.66\textwidth,width=1.06\textwidth]{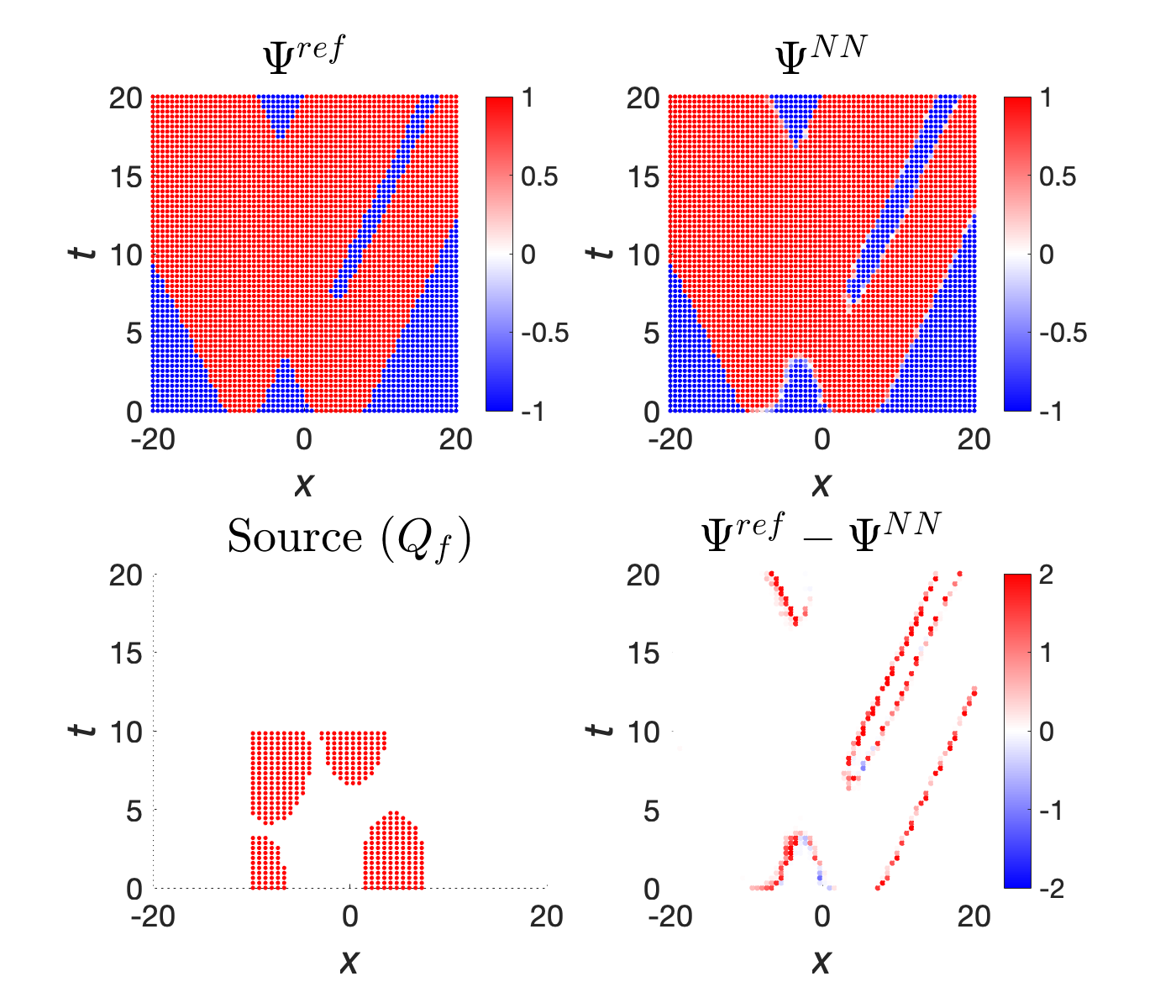}
    \includegraphics[trim=80 0 0 0,clip, height=0.66\textwidth,width=1.06\textwidth]{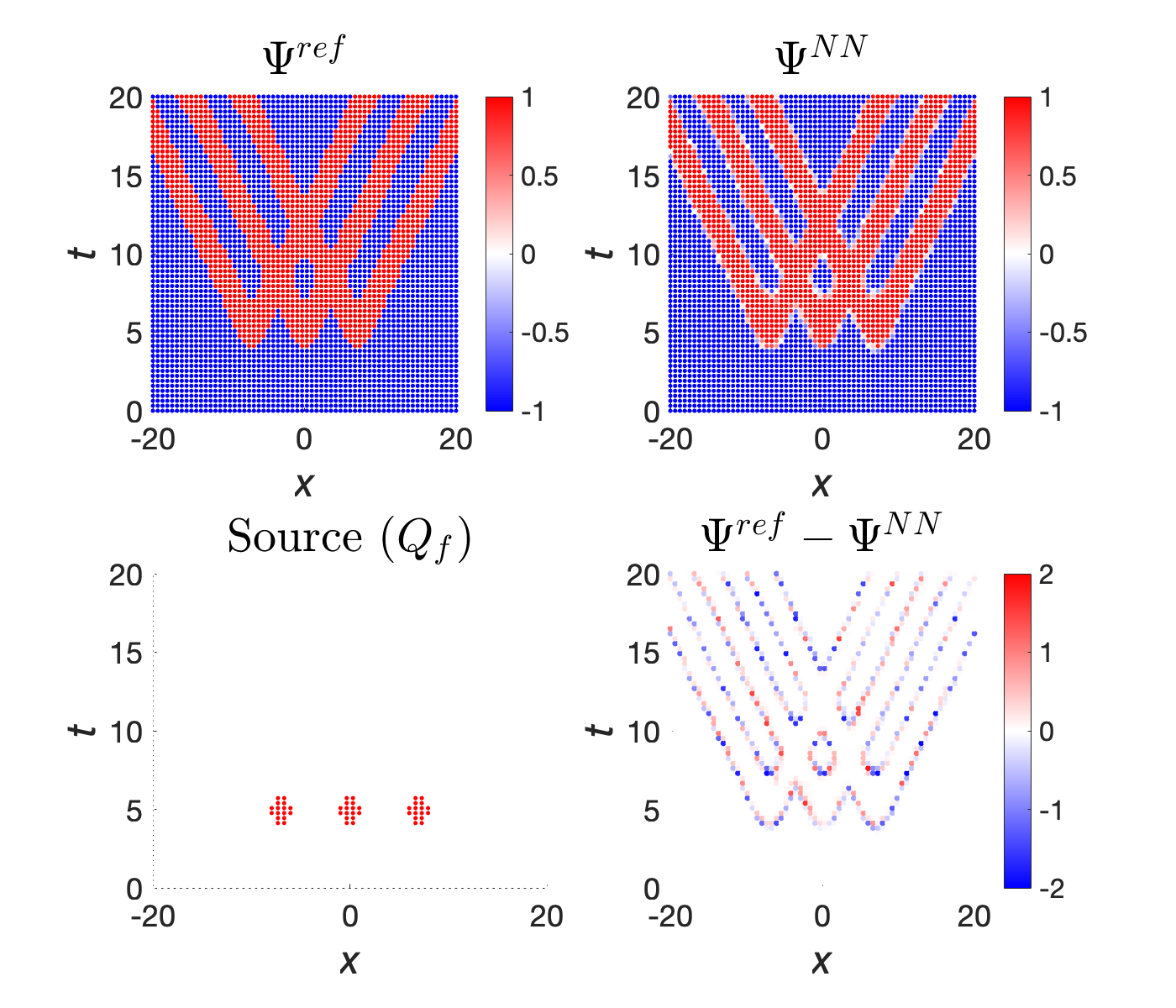}
%     \includegraphics[trim= 100 0 0 40,clip,width=1.6\textwidth]{./multi4}
%\end{center}
%\caption{Reconstruction of the shape of the lacuna (\ref{clt:lac1}) by neural network (\ref{clt:NN}) in the case $1\leq I^{(m)}\leq 4$ for two different examples. Example 1 is depicted in the top two rows and Example 2 is depicted in the bottom two rows with the same layout as in Figure \ref{clt:fig:single}.
\caption{Reconstruction of the shape of the lacuna (\ref{clt:lac1}) by neural network (\ref{clt:NN}) in the case $1\leq I^{(m)}\leq 4$ for two different examples. With the same layout as in Figure \ref{clt:fig:single}, the top two rows show an example from the test data set and the bottom two rows show a hand crafted example such that the lacunae has a 'pocket'.
\label{clt:fig:multi}}
\end{figure}

In this case, we generate our data choosing $I^{(m)}$ to be a random integer such that $1\leq I^{(m)}\leq4$ for each $m=1,\ldots,M$. Once trained, the neural network was able to predict which set, $\Lambda_1(Q_f)$ or $(\Lambda_1(Q_f))^c$, each node belongs to with approximately $98.55\%$ accuracy over $1{,}000$ test cases where the accuracy is calculated as in \eqref{clt:acc}. Figure \ref{clt:fig:multi} shows 2 examples (using the same layout as in Figure \ref{clt:fig:single}). The first example is taken from the test set, and in the second example we chose $Q_f$ such that the lacuna (\ref{clt:lac1}) has a ``pocket,'' i.e., a fully enclosed area. It is a part of the secondary lacuna (\ref{clt:lac1-2}). 

\section{Discussion} \label{clt:sec5}

We have demonstrated that a fully connected neural network can accurately reconstruct the shape of the lacunae in an artificial one-dimensional setting introduced in Section~\ref{clt:sec1}. While we have trained our network to find the shape of a combined lacuna 
(\ref{clt:lac1}), we anticipate that having it identify only the secondary lacunae (\ref{clt:lac1-2}) would not present any additional 
issues. A challenging next step is to extend the proposed machine learning approach to a realistic three-dimensional setting where the 
secondary lacunae are defined according to (\ref{clt:lac}) and account for the actual physics of the solutions to the wave equation 
(\ref{clt:wave}).

\begin{acknowledgement}
The work of A. Chertock and C. Leonard was supported in part by NSF Grant DMS-1818684. The work of S. Tsynkov was partially supported by 
the US-Israel Binational Science Foundation (BSF) under grant \#~2020128.
\end{acknowledgement}

\bibliographystyle{spphys}
\bibliography{references}
%\bibliography{../references}

\end{document}